\documentclass[10pt]{amsart}
\usepackage{amsfonts, amssymb,amsmath, amsthm, amsxtra, latexsym, amscd}

\newtheorem{theo+}    {Theorem}      [section]
\newtheorem{coro+}  [theo+]  {Corollary}
\newtheorem{lemm+}  [theo+]  {Lemma}

\theoremstyle{definition}
\newtheorem{ques+}   {Question}
\newtheorem{exam+}  [theo+]  {Example}
\newtheorem{rema+}  [theo+]  {Remark}

\newcommand\fg{\mathfrak g}
\newcommand\fk{\mathfrak k}

\newcommand\fp{\mathfrak p}
\newcommand\ad{\operatorname{ad}}
\let\Cal\mathcal
\newcommand\CC{{\mathbb C}}
\newcommand\DD{{\mathbb D}}

\def\[{\begin{equation}}
\def\]{\end{equation}}
\newcommand\HH{\Cal H}
\newcommand\RR{\mathbb R}
\newcommand\ZZ{\mathbb Z}
\def\poch(#1){(\!(#1)\!)}
\newcommand\LL{\mathcal L}
\newcommand\BB{\mathbb B}

\begin{document}

\title[Ramadanov conjecture]{Ramadanov conjecture and line bundles
over compact Hermitian symmetric spaces}

\author{Miroslav Engli\v s and Genkai Zhang}

\address{Mathematics Institute, Silesian University, 
Na~Rybn\'\i\v cku~1, 74601~Opava, Czech~Republic {\rm and} 
Mathematics Institute, Academy of Sciences, \v Zitn\'a 25, 
11567~Prague~1, Czech~Republic}

\email{englis@math.cas.cz}

\address{Department of Mathematical Sciences, Chalmers University 
of Technology, {\rm and} Department of Mathematical Sciences, 
G\"oteborg University, SE-412 96 G\"oteborg, Sweden}

\email{genkai@math.chalmers.se}

\thanks{Research supported by the Swedish Science Council (VR),
GA~\v CR grant 201/06/0128, and Czech Ministry of Education research 
plan MSM4781305904.}


\begin{abstract} We~compute the Szeg\"o kernels of the unit circle
bundles of homogeneous negative line bundles over a compact Hermitian
symmetric space. We~prove that their logarithmic terms vanish in all
cases and, further, that the circle bundles are not diffeomorphic to 
the unit sphere in $\CC^n$ for Grassmannian manifolds of higher ranks.
In~particular they provide an infinite family of smoothly bounded
strictly pseudo-convex domains on complex manifolds for which the log
terms in the Fefferman expansion of the Szeg\"o kernel vanish and which
are not diffeomorphic to the sphere. The~analogous results for the
Bergman kernel are also obtained.
\end{abstract}

\maketitle

\section{Introduction}
Let $\Omega$ be a strongly pseudo-convex bounded domain in $\CC^n$ with
smooth boundary. The~Bergman kernel has an expansion near the diagonal
in terms of the defining function of the domain, with the leading term
behaving like that of the Bergman kernel of the unit ball;
see~\cite{BFG-BAMS} and \cite{Monvel-Sj}. Similar result holds also for
the Szeg\"o kernel. However, there is in general also a logarithmic term
in the expansion of the Bergman and Szeg\"o kernel, and the study of the 
log term is of considerable interest for analytic and geometric motivations.
Among other things there is the Ramadanov conjecture \cite{Ram} which 
asserts that the answer to the following question is affirmative. 

\begin{ques+}\label{QUE1}
Let $\Omega$ be a strongly pseudo-convex bounded domain in $\CC^n$ with
smooth boundary. Suppose that the Bergman kernel has no logarithmic term. 
Is~the domain biholomorphic to the unit ball in $\CC^n$?
\end{ques+}

For certain special cases, such as domains in $\CC^2$, domains with
transversal or rotational symmetries, etc., this has been proved to
be~true; see~\cite{BCo}, \cite{Gra1}, \cite{Gra2}, \cite{BdM}, \cite{Nak}
and~\cite{Hir2}, and \cite{Att} for a real-variable version. 

There is also an obvious analogue of the conjecture for the Szeg\"o,
instead of Bergman, kernel, and one may also consider smoothly bounded
strictly pseudoconvex domains in complex manifolds. In~this setup,
in~particular, the following special case of the Ramadanov conjecture
was formulated in~\cite{Lu-Tian}.

\begin{ques+}\label{QUE2}
Let $S(L^\ast)$ be the disc bundle of a negative line bundle over a
simply-connected K\"ahler manifold~$M$. Suppose that the Szeg\"o kernel of
$S(L^*)$ has no log term. Is~the circle bundle diffeomorphic to the sphere?
\end{ques+}

However, as suggested to us by the anonymous referee,
since the singularity of the Bergman and Szego kernels are 
determined locally by the
CR structure of the boundary, a  reasonable related
conjecture
should be

\begin{ques+}\label{QUE3}
Let $\Omega$ be a strongly pseudo-convex  domain in a 
complex manifold with
smooth boundary. Suppose that the Bergman kernel has no logarithmic term. 
Is the boundary locally
CR equivalent to the sphere?
\end{ques+}

In~this paper, we~will consider the generating positive line bundle over a
compact Hermitian symmetric space $M$ and compute the corresponding Szeg\"o 
kernel. As~a~consequence we will see that the answer to the last question 
is negative. In~fact, the~simplest counterexamples are the powers $L^*=
\LL^{*m}$, $m>1$, of~the tautological bundle $\LL^*$ over the complex
projective space~$\CC P^n$: then the Szeg\"o kernel of $S(\LL^{*m})$ has 
no log term but $S(\LL^{*m})$ is the lens space $S^{2n+1}/\ZZ_m$ which is
not diffeomorphic to $S^{2n+1}$ for $m>1$. We~hasten to remark, however,
that these examples are not so interesting since they are still locally 
CR-equivalent to the sphere (i.e.~locally spherical). On~the other hand,
for~compact symmetric spaces of higher rank, we~even get examples 
which are not diffeomorphic to $S^{2n+1}/\ZZ_m$ for any $m\ge1$ --- in~fact,
they are not locally spherical at any point, providing counterexamples
to Question 3. Our~results thus indicate that 
some topological conditions are probably needed to have an affirmative 
answer to Question~\ref{QUE2}.%
\footnote{We~remark that there exist \emph{unbounded} pseudoconvex domains 
in $\CC^2$ with smooth boundary, as~well as bounded pseudoconvex domains 
in $\CC^2$ with \emph{rough} boundary, for which the Ramadanov conjecture 
for the Szeg\"o kernel is known to~fail; cf.~Remark~1.2 in~\cite{Hir1}.}

The~analogous assertions for the Bergman kernel are also established.
Here the simplest counterexample is the unit disc bundle $D(\LL^*)$
of the above-mentioned tautological bundle $\LL^*$ over~$\CC P^n$,
$n\ge1$, whose Bergman kernel has no log term but $D(\LL^*)$ is not
biholomorphic to the unit ball of $\CC^{n+1}$ 
(even though its boundary $S(\LL^*)$ is diffeomorphic to~$S^{2n+1}$).
Again, for~$\LL^*$ the tautological bundle over compact symmetric spaces 
of rank bigger than 1 we also get counterexamples whose boundary is not 
even locally spherical.

We~would like to thank Kengo Hirachi, Zhiqin Lu, Henrik Sepp\"anen, Robert
Stanton and Jan-Alve Svensson for several helpful and illuminating discussions.
The~authors also thank the referee for valuable comments.

\section{Compact Hermitian symmetric spaces}
We~briefly recall some necessary facts on compact Hermitian symmetric 
spaces; see~e.g.~\cite{He1}.

Let $\fg$ be a real simple Lie algebra of Hermitian type and let
$\fg=\fk+\fp$ be a Cartan decomposition of~$\fg$, where $\fk$ has
one-dimensional center $\RR Z$. Let
\[
\fg^{\CC} =\fp^- +\fk^{\CC} +\fp^+,   \label{2.1}
\]
where $\fp^{\pm}$ is the eigenspace of $\ad(Z)$ with eigenvalues $\pm i$.

Let $G^\CC$ be the simply connected Lie group with Lie algebra $\fg^\CC$, 
and let $G$, $K$, $P^{\pm}$ be the analytic subgroups of $G^\CC$ with Lie
algebras $\fg$, $\fk$, and $\fp^{\pm}$. Now $P^+K^\CC P^{-}$ is a dense
subset of $G^\CC$. For $g\in G^\CC, z\in \fp^{+}$ we let $g\cdot z$ and
$\Cal K(g: z)$ be the $\fp^+$ and $K^\CC$ components of $g\exp(z)$,
respectively. Namely,
\[
g\exp(z)=\exp(g\cdot z)\Cal K(g:z)p_{-}    \label{2.2}
\]
for some $p_{-}\in P^{-}$. Under the above action the $G$-orbit $G\cdot 0=G/K$ of $z=0\in \fp^{+}$ is a bounded domain in $\fp^{+}$, which
is the so-called Harish-Chandra realization of $G/K$.

Let $G^\ast$ be the analytic subgroup of $G^\CC$ with Lie algebra $\fg^\ast
= \fk+i\fp$. Then $M=G^\ast/K$ is a compact Hermitian symmetric space.
Furthermore $M=G^\CC/K^\CC P^{-}$, and under the identification of $\fp^+$
with $P^+$, the space $\fp^+$ is imbedded in $M$ as a dense subset. 
(In~fact, the complement of $\fp^+$ in $M$ is a complex submanifold 
of smaller dimension.)

Denote $\Cal K(g)=\Cal K(g:0)$. From (\ref{2.2}) we have 
\[
g=p_{+}\Cal K(g) p_{-},    \label{2.3}
\]
for $g\in P^+ K^\CC P^-$. 

For $z, w\in \fp^{+}$ we let $\Cal K(z, \bar w)$ be the $K^\CC$-part of
$\exp(-\bar w)\exp(z)$ as in (\ref{2.3}), namely $\Cal K(z, \bar w)=
\Cal K(\exp(-\bar w)\exp(z))$. The \emph{Bergman operator} $B(z, w)$ is
defined~by
$$
B(z, w)=\ad_{\fp^{+}}\Cal K(z, \bar w). 
$$
There exists an irreducible polynomial $h(z, w)$, called the 
\emph{Jordan canonical polynomial}, and an integer~$p$, the \emph{genus}
of~$G/K$ (see~(\ref{2.9}) below), such~that $\det B(z, w)=h(z, w)^{p}$; 
see~e.g.~\cite[\S4.15--4.17 and~\S7.4]{Loos-bsd}.

We~normalize the $K$-invariant inner product $(\cdot, \cdot)$ on $\fp^{+}$~by
$$
(z, w)=-\frac{1}{2p}\operatorname{tr} \left( \operatorname{Ad}_{\fp^+}(z) 
\operatorname{Ad}_{\fp^+}(\bar w) \right)
$$
where as before
\[
p=(r-1)a +2 +b     \label{2.9}
\]
is the {\it genus} of $M$. Note that the complex dimension of $M$ 
is $n=r+\frac{r(r-1)}2 a +rb$. We~let $\omega$ the $G^*$-invariant 
K\"ahler form on $M$ normalized so that the volume of $M$ with 
respect to $\omega^n$ is~$1$. In~terms of the local coordinates 
$z=\sum_j z_j e_j\in \fp^+$ (where $\{e_j\}$ is any 
orthonormal basis of $\fp^+$), viewing $\omega^n$ as a measure 
on $\fp^+$, we~have
\[
\omega(z)^n =C_0 \frac 1{h(z, -z)^p} \, dm(z),   \label{2.10}
\]
where $dm$ stands for the Lebesgue measure on~$\CC^n$. Further,
the~complement of $\fp^+$ in $M$ has zero measure with respect
to~$\omega^n$.

Finally we recall the Gindikin Gamma function and a version of
the generalized Pochhammer symbol defined~by
\begin{gather*}
\Gamma_M(c) = \prod_{j=1}^r\Gamma(c-\frac a 2(j-1)),   \\
\poch(c)_s=\frac{\Gamma_M(c+s)}{\Gamma_M(c)}, 
\end{gather*}
where $r$ is the rank of~$M$.

\section{Szeg\"o and Bergman kernels for the disc bundle}
Let~$L$ be the homogeneous line bundle over $M=G^*/K$ induced by 
the representation $k\mapsto (\det(\text{Ad} k))^{1/p}$, $k\in K$. 
(There exists a single-valued branch of the root of the determinant 
function so that this indeed defines a one-dimensional representation of~$K$.)
The bundle $L^p$ is then the top exterior product $\wedge^n {T^{(1, 0)}}$ 
of the holomorphic tangent bundle over $G^*/K$. Using the local coordinates 
$\fp^+\ni z\to \exp(iz)\in G^*/K$,
we~have that the fiber metric in $L^p$ is given~by
$$
\Vert\partial_1 \wedge \cdots \wedge \partial_n \Vert_{z}^2
= h(z, -z)^{-p}.
$$

Denoting by $e(z)$ a local holomorphic section of $L$ so that
$e(z)^p=\partial_1 \wedge \cdots \wedge \partial_n$ we see that 
the metric on $L$ is given by
$$
\Vert e(z) \Vert_{z}^2 = h(z, -z)^{-1}.
$$

Let $D=D(L^\ast)=\{\xi\in L^\ast; \Vert\xi\Vert < 1\}$ and
$S(L^\ast)=\partial D$ be the unit disc and the unit circle bundle 
of the dual bundle $L^*$ of~$L$, respectively. A~local defining 
function for $S(L^\ast)$ is given~by  
$$
\rho(z, \lambda e^\ast (z)) = {|\lambda|^2}{h(z, -z)} -1,
 \quad \lambda \in \CC, \; z\in\fp^+,
$$
where $e^\ast (z)$ is the local section of $L^\ast$ dual to $e(z)$. 
The~circle bundle $S(L^\ast)$ is a $CR$-manifold, with the $CR$-structure
defined by~$\rho$, and the disc bundle $D$ is strictly pseudoconvex,  
namely the Hessian $\partial\bar\partial\rho$ is positive definite 
on the holomorphic tangent space of~$S(L^\ast)$. 

The~manifold $S(L^\ast)$ is actually a compact homogeneous space of
$G^*\times S^1$, with $S^1=\{e^{i\theta}\}$ acting on $L^\ast$ fiberwise. 
Let $\pi$ be the projection $S(L^\ast)\to M$. We let $d\sigma$ be the
measure 
\[
d\sigma =\pi^\ast(\omega^n)\wedge \frac{d\theta}{2\pi}   \label{MEAS}
\]
which is also the unique $G^*\times S^1$-invariant probability measure
on~$S(L^\ast)$. Here $\omega $ is the K\"ahler form given above.

Let $\nu$ be a non-negative integer and consider the space $\HH^\nu$
of all holomorphic functions $\phi$ on $L^*$ satisfying
$$ \phi(e^{i\theta}\xi) = e^{\nu i \theta} \phi(\xi),
\qquad \xi\in L^*, \theta\in\RR.  $$
Note that owing to the holomorphy of $f$ this implies that even
$$ 
\phi(\lambda\xi) = \lambda^\nu \phi(\xi), \qquad \xi\in L^*, \lambda\in\CC.
$$ 
As~$M$ is compact, any such function~is, in~particular, automatically
square-integrable over $S(L^*)$ as well as over~$D$.

Identifying $\fp^+$ with a dense open subset of $M$ of full measure as
described in the previous section, and using the local trivializing
section $e^*(z)$ as before, the correspondence $\fp^+\times\CC\ni
(z,\lambda)\longleftrightarrow \xi=(z,\lambda e^*(z))\in L^*$ sets~up 
a bijection between a dense open subset of $D$ of full measure and 
the Hartogs domain  
$$  \Omega = \{ (z,\lambda)\in\fp^+\times\CC:
\; |\lambda|^2 h(z, -z) <1 \}.  $$
The~functions $\phi$ in $\HH^\nu$ then correspond to square-integrable
(with respect to~(\ref{2.10}) and the Lebesgue measure in~$\lambda$)
holomorphic functions $\tilde\phi$ on $\Omega$ satisfying $\tilde\phi
(z,\lambda)=\lambda^\nu f(z)$ for some entire function $f$ on~$\fp^+$.
\footnote{It~is clear that any function in $\HH^\nu$ must be of this 
form when restricted to the above local chart. Conversely, any square
integrable holomorphic function on $\Omega$ as above automatically extends
to a holomorphic function on all~of~$D$. Indeed, since the complement of
$\fp^+$ in~$M$ is a proper complex submanifold (see~e.g.~the discussion 
in \S2 in Berezin~\cite{Be}), making a suitable change of coordinates
it~is enough to show that any square-integrable holomorphic function
on the punctured polydisc $\DD^{n-1}\times(\DD\setminus\{0\})$ extends
to a holomorphic function on the whole~$\DD^n$. This ``$L^2$-version of
the removable singularity theorem'' is then easily proved by looking at
the Laurent expansion in~$z_n$, cf.~the proof for $n=1$ in~\cite{ACM}.}
\label{FOOTNOTE}
The~norm of $\phi$ in $L^2(d\sigma)$ thus equals~to
$$  \Vert f\Vert_{\nu}^2 = \int_{\fp^+} |f(z)|^2 h(z, -z)^{-\nu} 
\omega(z)^n .  $$
The~space $A^2_\nu(\fp^+)$ of all entire functions $f$ on $\fp^+$ for which
this norm is finite carries a representation of~$G^\ast$:
$$
g\in G^\ast: f(z)\mapsto f(g^{-1}z) J_{g^{-1}}(z)^{-\nu/p},
$$
where $J_{g^{-1}}$ is the complex Jacobian and $p$ is the genus
defined~in~(\ref{2.9}). 

The~function $h(z,-z)$ on~$\fp^+$ thus transforms according~to
\begin{align*}
h(g(z), -g(z))^{\nu} 
&= \frac{h(z,-z )^{\nu} h(w, -w)^{\nu} } {|h(z, -w)^{\nu}|^2}  \\
&= h(z,-z)^\nu \; |J_{g^{-1}}(z)|^{2\nu/p} ,  
\end{align*}
for $g\in G^\ast$ such that $g(0)=w$. 

\begin{lemm+} \label{LEMM}
The reproducing kernel for the space $A^2_\nu(\fp^+)$ is given~by
$$
\frac{\poch(\nu+p-\tfrac nr)_{\frac nr}} {\poch(p-\frac nr)_{\frac nr}}
\; h(z, -w)^\nu.
$$
\end{lemm+}

\begin{proof} It~follows from the transformation rule of $h(z, -w)$ 
under $G^\ast$ (see e.g.~\cite{Z}) that the reproducing kernel~is 
$$
c_{\nu} h(z, -w)^{\nu}.
$$
We~evaluate the constant, which is given~by the norm square of the
function~$1$, 
$$
c_{\nu}^{-1} = \int_{M} \Vert e^\nu\Vert ^2 \omega^n
= C_0 \int_{\fp^+} h(z, -z)^{-\nu }  h(z, -z)^{-p} \, dm(z) .
$$
In~terms of the polar coordinates (see \cite{FK}) we~have
$$
c_{\nu}^{-1} = C_0 C \int_{(\RR^+)^r} 
\prod_{j=1}^r (1+t_j^2)^{-\nu-p}    \prod_{j=1}^r t_j^{1+2b}
\prod_{1\le i<j\le r}^r |t_i^2-t_j^2|^a  \, dt_1\dots dt_r,
$$
with some constant $C$ independent of $\nu$. Changing variables~to
$t_j^2=x_j(1-x_j)^{-1}$, $j=1,\dots,r$, we~find that
$$
c_{\nu}^{-1} = \frac{C_0 C}{2^r} \int_{(0, 1)^r}
\prod_{j=1}^r (1-x_j)^{\nu}  \prod_{j=1}^r x_j^{b}
\prod_{1\le i<j\le r}^r |x_i-x_j|^a  \, dx_1 \dots dx_r ,
$$
which in turn can be expressed in terms of the Gindikin Gamma
function \cite{FK}, viz.,
$$
c_{\nu}^{-1} = C' \, \frac{\Gamma_M(\nu+p-\frac nr)} {\Gamma_M(\nu+p)}
= \frac {C'}{\poch(\nu +p- \tfrac nr)_{\frac nr}},
$$
with some constant $C'$ independent of~$\nu$. Taking $\nu=0$ and
recalling that $\omega^n$ was normalized to have total mass~one,
i.e.~$c_0=1$, gives $C'=\poch(p-\frac nr)_{\frac nr}$.
This completes the proof.
\end{proof}

Denote~by
$$
\rho(x, \alpha; y, \beta) = \alpha \bar \beta h(x, -y) - 1 
$$
the sesqui-holomorphic extension of the defining function $\rho$.

\begin{theo+}  \label{3.2}
The Szeg\"o kernel of the disc bundle $D$ is given, in local 
coordinates $\alpha e^*(z)\mapsto (z,\alpha) \in \fp^+\times\CC$, 
$|\alpha|^2 h(z,-z) < 1$,~by
\[
K(x, \alpha; y, \beta) = \sum_{\nu=0}^\infty 
\frac {\poch(\nu+p-\tfrac nr)_{\frac nr}} {\poch(p-\frac nr)_{\frac nr}}
\, h(x,-y)^{\nu} (\alpha\bar\beta)^{\nu}.    \label{Szego1}
\]
It~has an expansion in terms of the defining function $\rho$~as
\[
K(x, \alpha; y, \beta) = c_0 \rho(x, \alpha; y, \beta)^{-n-1} + 
c_1 \rho(x, \alpha; y, \beta)^{-n} + \cdots 
+ c_n \rho(x, \alpha; y, \beta)^{-1}    \label{Szego2}
\]
where $c_0=(-1)^{n+1}\dfrac{n!}{\poch(p-\frac nr)_{\frac nr}}$ and
$c_j$ are some real constants.
\end{theo+}

\begin{proof} It~is clear that the space $\HH^\nu$, $\nu=0,1,2,\dots$,
are~pairwise orthogonal subspaces of $L^2(d\sigma)$, and that their closed
span is the Hardy space $H^2(D)$. (In~fact, $H^2(D)=\oplus_{\nu=0}^\infty
\HH^\nu$ is just the Fourier decomposition of $H^2(D)$ into irreducible
components with respect to the action of~$S^1$.) Consequently, the Szeg\"o
kernel --- the~reproducing kernel of $H^2(D)$ --- is~the sum of the
reproducing kernels of the spaces $\HH^\nu$ over all~$\nu$, which
gives~(\ref{Szego1}). Since $\poch(\nu+p-\frac nr)_{\frac nr}$ is always 
a monic polynomial in $\nu$ of degree~$n$ --- hence, a~linear combination 
of the expressions $\frac{(\nu+1)\dots(\nu+k)}{k!}= \frac{(k+1)_\nu}
{\nu!}$, $k=0,1,2,\dots,n$ --- and
$$  \sum_{\nu=0}^\infty \frac{(k+1)_\nu}{\nu!} \, h(x,-y)^\nu \,
(\alpha\bar\beta)^\nu = [-\rho(x,\alpha;y,\beta)]^{-k-1},   $$
the formula (\ref{Szego2}) follows.   \end{proof}

\smallskip

Recall that the Bergman space of a complex manifold of dimension $n$ 
is in general defined as the space of all holomorphic $(n,0)$-forms 
$f$ such that  
\[
(-i)^{n} \int f \wedge \bar f < +\infty.   \label{Bg1}  
\]
The~Bergman kernel is then, by~definition, the~$(n,n)$-form
$$  \sum_m f_m \wedge \bar f_m ,  $$
where $\{f_m\}$ is any orthonormal basis of the Bergman space, with respect
to the inner product $\langle f,g\rangle$ obtained by replacing $\bar f$ in
(\ref{Bg1}) by~$\bar g$. The~sum is independent of the choice of the basis, 
etc.; see~e.g.~\cite{Koba}. Of~course, if~the manifold is just a domain
in~$\CC^n$, then by the identification
\[
f(z) \, dz_1\wedge\dots\wedge dz_n \longleftrightarrow f(z)   \label{Bg2}
\]
of $(n,0)$-forms with functions one recovers the usual definition of the
Bergman space and Bergman kernel of domains in~$\CC^n$.

We~have now a complete analogue of Theorem~\ref{3.2} also for the Bergman
kernel. 

\begin{theo+}  \label{3.3}
The Bergman kernel of the disc bundle $D$ is given, in local 
coordinates $\alpha e^*(z)\mapsto (z,\alpha) \in \fp^+\times\CC$, 
$|\alpha|^2 h(z,-z) < 1$,~by
$$  K(x, \alpha; y, \beta) = K^*(x,\alpha;y,\beta) \; dx_1 \wedge \dots 
\wedge dx_n\wedge d\alpha \wedge d\bar y_1 \wedge \dots \wedge 
d\bar y_n \wedge d\bar\beta,   $$
where 
\[
K^*(x,\alpha;y,\beta) = \frac1\pi \sum_{\nu=0}^\infty (\nu+1)
\frac{\poch(\nu+p-\tfrac nr)_{\frac nr}} {\poch(p-\frac nr)_{\frac nr}}
\, h(x,-y)^{\nu+1} (\alpha\bar\beta)^{\nu}.     \label{Bg3}
\]
It~has an expansion in terms of the defining function $\rho$~as
\[
\frac{K^*(x, \alpha; y, \beta)} {h(x,-y)}
= c_0 \rho(x, \alpha; y, \beta)^{-n-2} + 
\cdots + c_{n+1} \rho(x, \alpha; y, \beta)^{-1}    \label{Bg4}
\]
where $c_0=(-1)^{n+2}\dfrac1\pi\,\dfrac{(n+1)!} {\poch(p-\frac nr)
_{\frac nr}}$ and $c_j$ are some real constants.
\end{theo+}

\begin{proof} In~the local coordinates, we~can still identify the
$(n,0)$-forms with functions via~(\ref{Bg2}), and thus the Bergman
space on $D$ can be identified with the space of all functions
holomorphic and square-integrable on~$\Omega$, i.e.~with the usual
Bergman space on the Hartogs domain $\Omega\subset\CC^{n+1}$.~\footnote%
{Again, any~such $(n,0)$-form automatically extends to be holomorphic
even on the whole~$D$, i.e.~also on the complement of~$\Omega$ in~$D$;
see~the footnote before Lemma~\ref{LEMM} on~page~\pageref{FOOTNOTE}.}
Using again the Fourier decomposition with respect to the $S^1$-action,
together with the fact that now the norm of a function from~$\HH^\nu$,
$\nu=0,1,2,\dots$, equals
\begin{align*}
\|\phi\|^2_{L^2(\Omega)} &= \int_{\fp^+} \int_{|\lambda|^2<1/h(z,-z)}
|\lambda|^{2\nu} \, |f(z)|^2 \, d\lambda \wedge d\bar \lambda\,\wedge \omega(z)^n  \\
&= \frac\pi{\nu+1} \int_{\fp^+} |f(z)|^2 \, h(z,-z)^{-\nu-1}
\, \omega(z)^n,   \end{align*}
and, consequently, the reproducing kernel of $\HH^\nu$ with respect 
to this norm equals
$$ \frac{\nu+1}\pi \, \frac{\poch(\nu+1+p-\tfrac nr)_{\frac nr}}
{\poch(p-\frac nr)_{\frac nr}} \, h(z,-w)^{\nu+1}, $$
we~get the first formula (\ref{Bg3}) in the theorem. The~second formula
(\ref{Bg4}) follows from it in the same way as in Theorem~\ref{3.2}.
\end{proof}

With trivial modifications, the~last two theorems extend also to the 
unit circle bundles $S(L^{*\mu})$ and the corresponding unit disc bundles
$D_\mu=D(L^{*\mu})$ of the higher powers $L^{*\mu}$ of~$L^*$, $\mu=0,1,2,\dots$
(one~just needs to replace $\nu$ by $\nu\mu$ everywhere in the proofs.)

\begin{theo+}  \label{3.2A}
The Szeg\"o kernel of the disc bundle $D_\mu$ is given, in local 
coordinates $\alpha e^*(z)\mapsto (z,\alpha) \in \fp^+\times\CC$, 
$|\alpha|^2 h(z,-z)^\mu < 1$,~by
$$
K(x, \alpha; y, \beta) = \sum_{\nu=0}^\infty 
\frac {\poch(\mu\nu+p-\tfrac nr)_{\frac nr}} {\poch(p-\frac nr)_{\frac nr}}
\, h(x,-y)^{\mu\nu} (\alpha\bar\beta)^{\nu}. 
$$
It~has an expansion in terms of the sesqui-holomorphically extended 
defining function 
$\rho(x,\alpha;y,\beta)= \alpha\bar\beta h(x,-y)^\mu-1$~as
$$
K(x, \alpha; y, \beta) = c_0 \rho(x, \alpha; y, \beta)^{-n-1} + 
c_1 \rho(x, \alpha; y, \beta)^{-n} + \cdots 
+ c_n \rho(x, \alpha; y, \beta)^{-1} 
$$
where $c_0=(-1)^{n+1}\dfrac{n!\mu^n}{\poch(p-\frac nr)_{\frac nr}}$ and
$c_j$ are some real constants.
\end{theo+}

\begin{theo+}  \label{3.3A}
The Bergman kernel of the disc bundle $D_\mu$ is given, in local 
coordinates $\alpha e^*(z)\mapsto (z,\alpha) \in \fp^+\times\CC$, 
$|\alpha|^2 h(z,-z)^\mu < 1$,~by
$$  K(x, \alpha; y, \beta) = K^*(x,\alpha;y,\beta) \; dx_1 \wedge \dots 
\wedge dx_n\wedge d\alpha \wedge d\bar y_1 \wedge \dots \wedge 
d\bar y_n \wedge d\bar\beta,   $$
where 
$$
K^*(x,\alpha;y,\beta) = \frac1\pi \sum_{\nu=0}^\infty (\nu+1)
\frac{\poch(\mu\nu+p-\tfrac nr)_{\frac nr}} {\poch(p-\frac nr)_{\frac nr}}
\, h(x,-y)^{\mu\nu+1} (\alpha\bar\beta)^{\nu}. 
$$
It~has an expansion in terms of the sesqui-holomorphically extended 
defining function
$\rho(x,\alpha;y,\beta)= \alpha\bar\beta h(x,-y)^\mu-1$~as
\[
\frac{K^*(x, \alpha; y, \beta)} {h(x,-y)}
= c_0 \rho(x, \alpha; y, \beta)^{-n-2} + 
\cdots + c_{n+1} \rho(x, \alpha; y, \beta)^{-1}   \label{Bg4A}
\]
where $c_0=(-1)^{n+2}\dfrac1\pi\,\dfrac{(n+1)!\mu^n} {\poch(p-\frac nr)
_{\frac nr}}$ and $c_j$ are some real constants.
\end{theo+}

The~case of $\mu=p$ is of special interest, since in that case the 
$G^\ast\times S^1$-invariant probability measure (\ref{MEAS}) on 
$S(L^{*\mu})$ coincides with the surface measure used to get a 
holomorphically invariant Szeg\"o kernel, namely 
$$
\sigma\wedge d\rho = J[\rho]^{1/(n+2)} \, dV ,
$$
where $dV$ denotes the volume element in $\fp^+\times\CC$ and $J[\rho]$
stands for the Monge-Amp\'ere determinant
$$
J[\rho]=(-1)^{n+1} \det \bmatrix \rho & \partial\rho \\
\bar\partial\rho &\partial\bar\partial\rho \endbmatrix ;
$$
see e.g.~\cite{HK}. Indeed, a~short computation shows that $\sigma\wedge 
d\rho$ equals $h(z,-z)^{\mu-p}\,dV$ (up~to an immaterial constant factor),
while $J[\rho]=\mu^n h(z,-z)^{\mu-p}$; so~they coincide when $\mu=p$. 
Thus for $\mu=p$ Theorem~\ref{3.3A} concerns the invariant Szeg\"o kernel 
occurring in the theory of holomorphic invariants.

The~last four theorems yield abundant examples of smoothly bounded strictly
pseudoconvex domains in complex manifolds for which the Szeg\"o kernel as
well as the Bergman kernel contain no log-term in their boundary
singularity.  From the point of view of the Ramadanov conjecture,
it~remains to verify that these domains are not biholomorphic to the~ball.
Since by Fefferman's 1974 result \cite{Fef} any such biholomorphism
extends smoothly to the boundaries, it~is enough to show that the circle
bundle $S(L^{*\mu})$ is not diffeomorphic to the unit sphere~$S^{2n+1}$.

Recall that the simplest examples of compact Hermitian symmetric spaces
are the Grassmann manifolds $U(l)/U(k)\times U(l-k)$, $1\le k\le l-k$.
(They are the compact duals to the Cartan domains~$I_{k,l-k}$ --- the unit
balls $SU(k,l-k)/S(U(k)\times U(l-k))$ of complex $k\times(l-k)$ matrices.)
For $k=1$, the Grassmannians $M=U(l)/U(1)\times U(l-1)$
are just the complex projective spaces $M=\CC P^n$, $n=l-1$, and then 
the cosphere bundle $S(L^*)$ actually \underbar{is} $CR$-equivalent to 
the sphere $S^{2n+1}=U(l)/U(l-1)$: the bundle $L$ is the hyperplane 
bundle, $L^\ast$~is~the tautological bundle, and the mapping from 
the sphere $S^{2n+1}$ to $S(L^*)$ is given by $z\mapsto (\CC z, z)$.
Similarly, the~cosphere bundle $S(L^{*m})$ is $CR$-equivalent to the 
lens space $S^{2n+1}/\ZZ_m$, the~isomorphism now being given~by the mapping
$z\mapsto(\CC z,\otimes^m z)$ from the sphere $S^{2n+1}$ which induces 
a diffeomorphism from $S^{2n+1}/{\mathbb Z_m}$ onto $S({L^*}^m)$
(see~e.g.~\cite[p.~542]{Komu}).

From Theorems \ref{3.2}--\ref{3.3A} we thus arrive at the following
counterexamples to the manifold version of the Ramadanov conjecture
(Question~\ref{QUE2}).  

\begin{coro+}  \label{3.4a}
Let $M=\CC P^n$ be the complex projective $n$-space, $n\ge1$, and~$L$ the
positive line bundle as defined in the beginning of this section $($namely,
$L$~is the hyperplane bundle, i.e.~the dual of the tautological bundle$)$. 
Then the log-term vanishes in the Szeg\"o kernel of the circle
bundles~$S(L^{*\mu})$, and $S(L^{*\mu})$ is not diffeomorphic 
to the sphere $S^{2n+1}$ if $\mu>1$. 
\end{coro+}

\begin{proof} The~only thing we need to prove is that $S(L^{*\mu})\cong
S^{2n+1}/\ZZ_\mu$ is not diffeomorphic to $S^{2n+1}\equiv S^{2n+1}/\ZZ_1$ for
$\mu>1$. However, this is immediate for instance from the cohomology groups
(see e.g.~\cite[Example~18.5]{Bott-Tu} or \cite[Example 2.43, p.~144]{Hatcher})
\[ H^j(S^{2n+1}/\ZZ_\mu,\ZZ) = \begin{cases}  \ZZ, & j=0,2n+1, \\  \ZZ_\mu, 
& j=2, 4,\dots,2n, \\ 0 & \text{otherwise},  \end{cases}  \label{Hj}  \]
since the cohomology rings are diffeomorphic invariants.  \end{proof}

It~is not difficult to see that for Grassmannians of higher rank, we~even get
counterexamples which are not diffeomorphic to any lens space $S^{2n+1}/\ZZ_m$.

\begin{coro+}  \label{3.4b}
Let $M=U(l)/U(k)\times U(l-k)$ $(1<k \le l-k)$ be the Grassmannian 
of higher rank $k>1$ and complex dimension $n=k(l-k)$. Let $L$ be 
the positive line bundle as defined in the beginning of this section. 
$($Namely, $L$~is the determinant bundle of the hyperplane bundle$.)$ 
Then the log-term vanishes in the Szeg\"o kernel of the circle
bundles~$S(L^{*\mu})$, $\mu\ge1$, and $S(L^{*\mu})$ is not 
diffeomorphic to any lens space $S^{2n+1}/\ZZ_m$. 
\end{coro+}

\begin{proof} We~use the Gysin exact sequence 
\cite[Theorem 12.2]{Milnor-Stasheff} \cite[p.~437~ff.]{Hatcher} 
of the circle bundle $E:=S(L^{*\mu})$ over~$M$
(all~cohomology groups are over $\mathbb R$): 
$$
\cdots \to H^{2j-1}(E) \to H^{2j-2}(M) \to H^{2j}(M) \to H^{2j}(E) \to \cdots.
$$
If~$E$~were diffeomorphic to $S^{2n+1}/\ZZ_m$, then by (\ref{Hj}) we would have
$$ H^j(E) = \begin{cases} \RR, & j=0,2n+1,  \\ 0 & \text{otherwise}.
\end{cases}  $$
From the Gysin sequence it would thus follow that
\[  H^{2j-2}(M) \cong H^{2j}(M), \qquad j=1,\dots,n.    \label{HH}    \]
On~the other hand, it~is known that the Poincar\'e series of the 
cohomology ring $H^*(M)$ is given~by (see e.g.~\cite[Chapter~IV, 
Proposition~23.1]{Bott-Tu})
$$
\frac{(1-t^2)\cdots (1-t^{2l})}
{(1-t^2)\cdots (1-t^{2k}) (1-t^2)\cdots (1-t^{2(l-k)}) }.
$$
Thus (\ref{HH}) can happen only for $k=1$.   \end{proof}

The~lowest-dimensional counterexample to the Ramadanov conjecture for the
Szeg\"o kernel of circle bundles, namely Question~\ref{QUE2}, supplied by  
Corollary~\ref{3.4a} thus occurs for the circle bundles $S(L^{*m})$, $m>1$,
of powers of the tautological bundle over the Gauss sphere~$\CC P^1$ (so~that
$S(L^{*m})$ has real dimension~3), while that supplied by Corollary~\ref{3.4b}
--- i.e.~not diffeomorphic to the lens spaces ---  for~the Grassmannian with
$k=l-k=2$ (i.e.~with $S(L^*)$ of real dimension~9). 

Finally, we~also have the corresponding assertions for the Bergman,
instead of the Szeg\"o, kernel.

\begin{coro+}  \label{3.4c}
Let $M=U(l)/U(k)\times U(l-k)$ $(1\le k \le l-k)$ be the Grassmannian of rank
$k\ge1$ and complex dimension $n=k(l-k)$. Let $L$ be the positive line bundle
as defined in the beginning of this section. Then the log-term vanishes in 
the Bergman kernel of the corresponding disc bundles~$D_\mu$, $\mu\ge1$, and
$D_\mu$ is not biholomorphic to the unit ball $\BB^{n+1}$ of~$\CC^{n+1}$.  
\end{coro+}

\begin{proof} That $D_\mu$ is not biholomorphic to $\BB^{n+1}$ if $k>1$ or
$\mu>1$ follows from the last two corollaries since its boundary $\partial
D_\mu=S(L^{*\mu})$ is then not diffeomorphic to $\partial\BB^{n+1}=S^{2n+1}$.
We~claim that $D_\mu$ is still not biholomorphic to $\BB^{n+1}$ even if
$k=\mu=1$, i.e.~for $L^*$ the tautological line bundle over $M=\CC P^n$ 
(even though $S(L^*)$ then is diffeomorphic to~$S^{2n+1}$).
Indeed, a~short computation using (\ref{Bg4A}) shows that the zero 
section of $D(L^*)$ is then a totally geodesic submanifold with respect 
to the Bergman metric; since any biholomorphism is automatically 
an isometry with respect to Bergman metrics, it~would follow that
the image of the zero section is a compact submanifold of the unit
ball which is totally geodesic with respect to the Bergman metric.
However, no such submanifold can exist, since every geodesic in 
the ball with respect to the Bergman metric reaches the boundary
(the~geodesics through the origin are just straight lines, and the 
ball is homogeneous). Thus~$D(L^*)$ cannot be biholomorphic to the ball.
(This is in apparent contrast with the situation for domains in~$\CC^n$, 
where by the recent theorem of Chern and~Ji~\cite{CJ}, any smoothly-bounded
simply connected domain whose boundary is locally spherical must be
biholomorphic to the~ball.)
This completes the proof.   \end{proof}

In~particular, for $n=\mu=1$ the disk bundle $D$ over the Gauss sphere 
$\CC P^1$ provides a~two-dimensional counterexample to the manifold version 
of the Ramadanov conjecture for the Bergman kernel (Question~\ref{QUE1}). 

As~noted in the Introduction, the simplest counterexamples mentioned above
(i.e.~the circle bundles over $\CC P^n$ and the disc bundle over~$\CC P^1$)
are not so interesting, since it is apparent that $S(L^{*\mu})\cong S^{2n+1}/
\ZZ_\mu$ is still locally spherical (i.e.~locally CR-equivalent to the sphere).
In~fact, using the biholomorphism $z\mapsto(\CC z,z)$ already mentioned, it~is
seen that $D(L^*)$ with zero section removed is then biholomorphic to $\BB^n
\setminus\{0\}$, i.e.~$D(L^*)$ is just the one-point blow-up of the ball.
Remarkably, for compact symmetric space $M$ of higher rank the situation 
is already different.   

\begin{coro+} \label{BuSh} 
Let~$M$ and $L$ be as in Corollary~\ref{3.4b}. Then the log-term vanishes 
in the Bergman kernel of the corresponding disc bundles~$D_\mu$, $\mu\ge1$,
as~well as in the Szeg\"o kernel of the circle bundles~$S(L^{*\mu})$, and at 
the same time the boundary $S(L^{*\mu})$ of $D_\mu$ is not locally spherical
at any point.    \end{coro+}

\begin{proof} By~homogeneity, if~$S(L^{*\mu})$ were locally spherical at
some point, then it would be such at all points, i.e.~would be spherical.
Since $S(L^{*\mu})$ is compact and homogeneous, Proposition~5.1 in Burns 
and Shneider~\cite{BS} would then imply that $S(L^{*\mu})$ is isomorphic 
to the lens space $S^{2n+1}/\ZZ_m$ for some~$m$. But~we have seen in
Corollary~\ref{3.4b} that this is not the~case.   \end{proof}

\begin{rema+} We~remark that the disc bundles $D_\mu$ are not Stein.
(In~fact, any~holomorphic function on $D_\mu$ must necessarily be 
constant on the compact manifold~$M$, embedded as the zero section 
in~$D_\mu$.) The~authors do not know any counterexample to the 
conjecture which would be Stein.   \end{rema+}

\smallskip

In~view of the above results, it~seems somewhat natural to pose the
following modified version of the Ramadanov conjecture.

\begin{ques+}
Suppose that the Szeg\"o or Bergman kernel of a domain in a complex
manifold has no log term in its boundary singularity. Is~the domain
then always biholomorphic to the unit disc bundle $D(L^*)$ for some 
positive line bundle $L$ over a compact Hermitian symmetric space~$M$?
\end{ques+}

We~conclude by remarking that there is a well-known intimate relationship
between \emph{functions} on the dual disc bundle $D\subset L^*$ and
\emph{sections} of the tensor powers $L^\nu$ of the original line
bundle~$L$. Namely, let~$\Cal L^\nu$, $\nu=0,1,2,\dots$, stand for 
the space of all functions $f$ on $L^*$ satisfying 
$$ f(e^{i\theta}\xi) = e^{\nu i\theta} f(\xi),
 \qquad \xi\in L^*, \theta\in\RR. $$
Then the natural mapping $s\mapsto\tilde s$,
$$  \tilde s(\xi) := \langle s, \xi^{\otimes\nu} \rangle, 
\qquad \xi\in L^*,   $$
sets~up a bijection between functions $\tilde s\in\Cal L^\nu$ and sections
$s$ of~$L^\nu$; further, $s$~is holomorphic if and only if $\tilde s$~is
(i.e.~if~and only if $\tilde s$ belongs to the space~$\HH^\nu$).
In~this~way, some of the results in this paper can be recast in the
language of reproducing kernels of Bergman spaces of sections of the
powers $L^\nu$ of the line bundle~$L$. We~omit the details.


\end{document}